\documentclass[letterpaper, 10 pt, conference]{ieeeconf}  

\IEEEoverridecommandlockouts                              
\usepackage{graphics} 
\usepackage{epsfig} 
\usepackage{amsmath} 
\usepackage{amssymb}  
\usepackage{float}
\usepackage{graphicx}
\usepackage{caption}
\usepackage{subfig}
\usepackage[T1]{fontenc}
 
\newtheorem{lemma}{Lemma}
\newtheorem{theorem}{Theorem}

\newtheorem{definition}{Definition}

\newtheorem{remark}{Remark}

\usepackage{mathtools}
\usepackage{soul}
\usepackage{tikz}
\usepackage{color}
\usepackage{tcolorbox}
\mathtoolsset{showonlyrefs}

\title{\LARGE \bf
Controllability and Stabilizability Analysis \\ of Signed Consensus Networks
}

\author{Siavash Alemzadeh\thanks{The research of the authors has been supported by NSF grant SES-1541025 and AFOSR grant FA9550-16-1-0022. Siavash Alemzadeh is with the Mechanical Engineering and Mathias Hudoba de Badyn and Mehran Mesbahi are with the Department of Aeronautics and Astronautics, University of Washington, WA 98195. Emails: \{alems,hudomath,mesbahi\}@uw.edu. \textcopyright~IEEE 2017}, Mathias Hudoba de Badyn, and Mehran Mesbahi }%

\begin{document}
	
\maketitle
\thispagestyle{empty}
\pagestyle{empty}

\begin{abstract}

Signed networks have been a topic of recent interest in the network control community as they allow studying antagonistic interactions in multi-agent systems.
Although dynamical characteristics of signed networks have been well-studied, notions such as controllability and stabilizability for signed networks for protocols such as consensus are missing in the literature.
Classically,  graph automorphisms with respect to the input nodes have been used to characterize uncontrollability of consensus networks.
In this paper, we show that in addition to the graph symmetry, the topological property of structural balance facilitates the derivation of analogous sufficient conditions for uncontrollability for signed networks.
In particular, we provide an analysis which shows that a gauge transformation induced by structural balance allows symmetry arguments to hold for signed consensus networks.
Lastly, we use fractional automorphisms to extend our observations to output controllability and stabilizability of signed networks.

\vspace{2mm}

\textit{Index Terms}$-$Consensus networks,  controllability, stabilizability, signed
graphs, structural balance, graph symmetry

\end{abstract}

\section{INTRODUCTION}

Networked systems have been at the forefront of active research in the systems and control
community for the past two decades.
Some well-studied examples of networked systems are social networks and dynamics of opinions~\cite{Hegselmann2002}, flocking~\cite{vicsek1995}, autonomous robotics~\cite{Stone2000}, quantum networks~\cite{Acin2007}, autonomous flight~\cite{Model2000}, traffic control \cite{Yang1995} and gene regulatory networks~\cite{Yeung2002}.
Several networked systems include both cooperative and antagonistic interactions, such as certain classes of social dynamics~\cite{cartwright1956structural}.

Consensus algorithms have been used in many scientific and engineering applications, including multi-agent systems~\cite{Chen2013a,Saber2003}, robotics~\cite{Joordens2009} and Kalman filtering~\cite{Olfati-Saber2005}.
A large amount of research has been dedicated to looking at the control of consensus~\cite{Mesbahi2010}.
The work by Rahmani \emph{et al.} showed that certain symmetries of networks characterized by automorphisms of the topology of the network cause uncontrollability~\cite{Rahmani2009a}.
This was generalized by Chapman and Mesbahi who showed signed fractional automorphisms generate necessary and sufficient conditions for uncontrollability and unstabilizability of linear networks~\cite{Chapman2015a}.
Further works in this direction have examined methods of generating network topologies that are controllable for consensus, such as in~\cite{Chapman2014a} and~\cite{Hudobadebadyn2016}, or creating networks to satisfy performance measures \cite{Siami2017}, \cite{HudobaDeBadyn2015a}.

Consensus algorithms on networks with antagonistic interactions were first considered by~Altafini~\cite{Altafini2013,Altafini2012}.
The network property of \emph{structural balance}, first considered in the study of social networks (\cite{cartwright1956structural,Akiyama1981,Kaszkurewicz2012})
was identified in Altafini's work as the property inducing bipartite consensus in which the agents converge to two disjoint clusters instead of a uniform consensus.
Graph-theoretic properties of signed Laplacian dynamics were studied by Pan \emph{et al.}~\cite{Pan2016a}.
Further research by Pan \emph{et al.} has looked at identifying the bipartite structure of structurally balanced graphs using data from signed Laplacian dynamics and dynamic mode decomposition~\cite{Pan2016}, adding to the works done by Harary and Kabell~\cite{Harary1980} and Facchetti \emph{et al}~\cite{Facchetti2011}.
Recent contributions by Clark \emph{et al.} have studied the leader selection problem in signed consensus~\cite{Clark2017}.

The contributions of this paper are as follows.
We first conduct a controllability analysis of signed Laplacian consensus using symmetry arguments developed by Rahmani \emph{et al.} \cite{Rahmani2009a}, for the single-input-single-output (SISO) and multiple-input-multiple-output (MIMO) cases of consensus dynamics with leader nodes.
In particular, we identify the property of structural balance, that when combined with network symmetry,
causes uncontrollability of signed consensus dynamics.
The key feature of structurally balanced graphs that allows this analysis is that they admit a \emph{gauge transformation} that allows the permutation matrix corresponding to the graph symmetry be extended to a signed permutation matrix.
We then use tools developed by Chapman and Mesbahi in~\cite{Chapman2015a} to derive controllability and stablilizability conditions for influenced signed consensus dynamics.

The paper is organized as follows.
In \S\ref{sec:math-prel}, we discuss the mathematical background needed in the paper.
We outline the problem statement in \S\ref{sec:problem-statement}.
In \S\ref{sec:methods}, we show that structural balance combined with symmetry about inputs leads to uncontrollability, 
and derive the corresponding stabilizability conditions.
Relevant examples are shown in \S\ref{sec:examples}, and the paper is concluded in \S\ref{sec:disc-concl}.


\section{MATHEMATICAL BACKGROUND}

\label{sec:math-prel}

We consider $\mathbb{R}_+$ as the set of nonnegative real numbers.
A column vector with $n$ elements is designated as $v\in\mathbb{R}^n$ where $v_i$ or $[v]_i$ both represent its $i$th element.
Matrix $M\in \mathbb{R}^{p\times q}$ contains $p$ rows and $q$ columns with $[M]_{ij}$ denoting the $i$th row and $j$th column element of $M$.
For $w\in\mathbb{R}^n$ the $\mathrm{diag}(w)$ is an $n\times n$ matrix with $w$ on its diagonal and zero elsewhere.
The unit vector $e_i$ is the column vector with all zero entries except $[e_i]_i=1$.
The column vector of all ones is denoted as $\textbf{1}$.
The \emph{cardinality} of a set $S$ is denoted as $|S|$.
We define $S^{\perp}=\{ v^*\in \mathbb{R}^n : \langle v,v^* \rangle=0 \ \text{for all} \ v\in S \}$ 
where $\langle \cdot , \cdot \rangle$ is the inner product in the Euclidean space.
The column space of a matrix $M$ is denoted by $\mathcal{R}(M)$.
We define $\mathcal{R}(P)$ to be $A$-\emph{invariant} if there exists $C$ such that $AP=PC$. We say $A$ is similar to $B$ if there is an invertible matrix $R$ such that $R^{-1}AR=B$. Two similar matrices share the same spectra. The \emph{leading principal submatrix} of order $k$ of $X$ is the square submatrix of $X$
formed by deleting the last $n-k$ rows and columns.

\subsection{Consensus Dynamics in Signed Networks}\label{sec:graphss-graphsl}

A multi-agent system with $n$ agents is characterized by a \emph{graph} $\mathcal{G}=(\mathcal{V},\mathcal{E},{W})$ where $\mathcal{V}=\{v_1,v_2,\dots,v_n\}$ is the set of nodes, $\mathcal{E}\subseteq \mathcal{V}\times\mathcal{V}$ denotes the set of edges, and ${W}\in\mathbb{R}_+^{n\times n}$ consists of weights assigned to edges.
We say $j\in\mathcal{N}(i)$ if $i$ and $j$ are neighbors, i.e. there exists an edge between $i$ and $j$.
A \emph{path} of length $r$ in $\mathcal{G}$ is given by a sequence of different nodes $v_{i_0},v_{i_1},\dots,v_{i_r}$ such that $v_{i_k}$ and $v_{i_{k+1}}$ are neighbors for $k=0,1,\dots,r-1$.
When the terminal nodes are equal, the path is called a \emph{cycle}.
Also, $\mathcal{G}$ is \emph{connected} if for every $i$ and $j$, there is a path between $v_i$ and $v_j$.
The square matrix ${A}\in\mathbb{R}_+^{n\times n}$ is called the \emph{adjacency matrix} if ${A}_{ij}={W}_{ij}\neq 0$.
The \emph{degree matrix} ${D}\in\mathbb{R}^{n\times n}$ is a square diagonal matrix where ${D}_{ii}=\sum_{j\in\mathcal{N}(i)} {A}_{ij}$.
The \emph{graph Laplacian} is then defined as ${L}={D}-{A}$ and the \emph{consensus dynamics} $\dot{x}=-{L}x$~\cite{Mesbahi2010}.

A \emph{signed graph} $\mathcal{G}_s$ is a graph that admits negative weights.
The \emph{signed graph Laplacian} is defined as  ${L}_s={D}_s-{A}_s$, where ${A}_s$ can also contain negative elements and the degree matrix is again diagonal with $[{D_s}]_{ii}=\sum_{j\in\mathcal{N}(i)} \left|[{A}_s]_{ij}\right|$.
The corresponding dynamical system is given by $\dot x = -{L}_sx$, i.e. $\dot x_i = -\sum_{j\in\mathcal{N}(i)} \left|{W}_{ij}\right|\left(x_i - \mathrm{sgn}({W}_{ij})x_j\right)$
where $\mathrm{sgn}$ represents the sign function.
A \emph{positive cycle} is a cycle with even number of negative edges.
A \emph{gauge transformation} is a change of orthant order via a square matrix $G_t \in \left\{ \mathrm{diag}(\sigma):~\sigma=[\sigma_1,\dots,\sigma_n]~,~\sigma_i = \pm 1\right\}$.
Then $G_t=G_t^T=G_t^{-1}$ \cite{Altafini2012}.

\subsection{Automorphisms/Interlacing/Equitable Partitions}\label{sec:automorphisms}

An \emph{automorphism} of the graph $\mathcal{G}$ is a permutation $\psi$ of its nodes such that $\psi(i)\psi(j)\in\mathcal{E}$ if and only if $ij\in \mathcal{E}$.
Let the permutation matrix $\Psi$ be such that $[\Psi]_{ij}=1$ if $\psi(i)=j$ and zero otherwise.
Then $\psi$ is an automorphism of $\mathcal{G}$ if and only if $\Psi A(\mathcal{G})=A(\mathcal{G}) \Psi$ (see \cite{Mesbahi2010}).

Suppose $A\in\mathbb{R}^{n\times n}$ and $B\in\mathbb{R}^{m\times m}$ are both symmetric and $m\leq n$.
Then the eigenvalues of $B$ \emph{interlace} the eigenvalues of $A$ if for $i=1,2,\dots,m$, $\lambda_{n-m+i}(A)\leq \lambda_i(B)\leq\lambda_i(A)$ where $\lambda_1(A)\geq\lambda_2(A)\geq\dots\geq \lambda_n(A)$ are the eigenvalues of $A$ in a non-increasing order \cite{godsil2013algebraic}.

The \emph{cell} $C$ is a subset of the graph nodes $\mathcal{V}$. A \emph{nontrivial cell} is a cell with more than one node.
A \emph{partition} is a grouping of $\mathcal{V}$ into different cells.
An \emph{r-partition} $\pi$ of ${\mathcal{V}}$ with cells $\{C_i\}_{i=1}^r$ is \emph{equitable} if each node in $C_j$ has the same number of neighbors in $C_i$, for all $i,j$.
We call $\pi$ a \emph{nontrivial equitable partition} (NEP) if it contains at least one nontrivial cell.
Let $b_{ij}$ be the number of neighbors in $C_j$ of a node in $C_i$. 
The \emph{quotient} of $\mathcal{G}$ over $\pi$, denoted by $\mathcal{G}/\pi$, is the directed graph with the cells of an equitable \emph{r-partition} $\pi$ as its nodes and $b_{ij}$ edges directed from $C_i$ to $C_j$.
The adjacency matrix of the quotient is specified by $[A(\mathcal{G}/\pi)]_{ij}=b_{ij}$. A \emph{characteristic vector} $p_i\in \mathbb{R}^n$ of a nontrivial cell $C_i$ has 1's in components associated with $C_i$ and 0's elsewhere.
A \emph{characteristic matrix} $P\in \mathbb{R}^{n\times r}$ of a partition $\pi$ of $\mathcal{V}$ is defined as $[p_i]_{i=1}^r$.
More details and examples can be found in \cite{Rahmani2009a}.


\section{PROBLEM STATEMENT}\label{sec:problem-statement}

\label{sec:problem-statement}

The analysis of this paper consists of two main parts, which examine two different notions of control on networks.
In the first part, we examine the notion of uncontrollability which was initially derived in \cite{Rahmani2009a} but for signed consensus networks.
The notion of control in this case is taking over the state of one or several nodes in the network, and using their edges to inject signals into the system.
In the second part, we consider the case where the nodes are controlled by injecting a single-integrator signal to some nodes of the graph.

Signed consensus networks are of interest because the negative weights induce a phenomenon known as \emph{clustering}, where agents will not converge to an agreement subspace, but rather converge to opposite equilibria.
The condition that causes clustering was identified by~\cite{Altafini2013} as \emph{structural balance}.
We will show that this topological feature is the additional condition that produces uncontrollability of signed networks.

One important result from \cite{Altafini2013} that is frequently utilized is the following equivalences:
\begin{enumerate}
	\item The signed graph $\mathcal{G}$ is structurally balanced;
	\item There exists a gauge transformation $G_t$ such that $G_tA_sG_t$ has positive entries (i.e. $G_tA_sG_t$ is \emph{unsigned});
	\item All cycles in $\mathcal{G}$ are positive;
	\item The signed Laplacian $L_s$ has a zero eigenvalue;
	\item There exists a bipartition of $\mathcal{V}$ such that the edges within the same set are positive, and the connecting edges are negative.
\end{enumerate}

\noindent It is shown that if the signed graph is structurally unbalanced, $\lim_{t\to\infty} x(t) = 0$. Otherwise, $\lim_{t\to\infty}x(t) = (1/n) \left( \textbf{1}^T G_tx(0)\right) G_t \textbf{1}$ implying the bipartition.

\vspace{1mm}

Previous work by Rahmani \emph{et al.} in~\cite{Rahmani2009a} showed that symmetry with respect to a single input and interlacing for multiple input are sufficient for uncontrollability, and this was generalized by Chapman and Mesbahi in~\cite{Chapman2015a} to show that fractional symmetry with respect to the inputs is sufficient and necessary for uncontrollability. In this paper, we show that structural balance is the property that combined with symmetry and interlacing leads to uncontrollability.

Before proceeding to our main results, we summarize the various dynamics considered in the paper.

Given a connected signed graph $\mathcal{G}_s$, we can select one node and use it to inject our input signal $u$.
This corresponds to partitioning the Laplacian as follows:
\begin{align}
L_s = \left[ 
\begin{array}{c|c}
A_s^f&B_s^f\\\hline
{B_s^f}^T&A_s^i
\end{array}\right],
\end{align}
where $f$ and $i$ denote the \emph{floating} and \emph{input} parts of the network respectively.
Then the dynamical system for signed consensus networks is
\begin{align}
	\dot x = -A_s^f x -B_s^fu \ ,
	\label{eq:1}
\end{align}
which holds for both SISO and MIMO systems and the \emph{floating signed graph} is also denoted by $\mathcal{G}_s^f$.
Then we use the \emph{Popov-Belevitch-Hautus} (PBH) controllability test which specifies that the system is controllable if and only if none of the eigenvectors of $A_s^f$ are simultaneously orthogonal to all columns of $B_s^f$.
More details on this can be found in \cite{Mesbahi2010}.

The second variant of controlling consensus networks is to simply inject signals into nodes, without taking over the state of the node.
We can therefore define the \emph{influenced signed consensus dynamics} with $q$ inputs and $p$ outputs as
\begin{align}
	\dot x = -L_s x + B(I)u,~y = C(O)x,
	\label{eq:2}
\end{align}
where $I\subseteq N$ is the set of input nodes, $O\subseteq N$ the set of output nodes, and $B(I),C(I)$ are the matrices $B(I) = \begin{bmatrix}
e_{i_1}&\cdots&e_{i_q}
\end{bmatrix}
$ and $C(O) = \begin{bmatrix}
f_{j_1}&\cdots&f_{j_p}
\end{bmatrix}^T$
in which $e_i$ is the unit vector for nodes $i\in I$, and $f_j$ for nodes $j\in O$.
The control signal vector is $u\in\mathbb{R}^q$.
Both (output) controllability and (output) stabilizability of these dynamics is considered in Section~\ref{sec:stab-outp-contr}.


\section{Analysis}\label{sec:methods}

\subsection{Leader-Follower System Controllability}

\vspace{1mm}

\noindent\subsubsection{The SISO Case}\label{sec:sign-graph-contr}

~

\vspace{1mm}

A result by \cite{Rahmani2009a} shows that for a SISO consensus network, a symmetry about the input node is sufficient for uncontrollability.
We extend this result for the signed consensus networks considered by \cite{Altafini2013} and \cite{Pan2016}.
In particular, we show that structural balance and input symmetry is sufficient for uncontrollability of signed graphs.

\vspace{2mm}

\begin{remark}
	\label{rmk:1}
	This result shows that signed Laplacian is in some sense more robust to symmetries about the input nodes.
	In particular, there are examples of unsigned graphs that are symmetric about an input and therefore uncontrollable, but whose signed counterparts exhibit controllability despite the symmetry.
	Therefore, input symmetry does not impose enough constraints on the network to claim its uncontrollability.
\end{remark}

\vspace{2mm}

\begin{lemma}
	\label{lem:2}
	Assume the unsigned graph $\mathcal{G}$ enjoys input symmetry and the signed network $\mathcal{G}_s$ is structurally balanced. Then there exists a nontrivial $J'$ such that (1) $J'A_s^f=A_s^fJ'$, (2) $J'^TB_s^f=B_s^f$, and (3) if $v$ is the eigenvector corresponding to the eigenvalue $\lambda$ of $A_s^f$, then $J'v$ (and hence $v-J'v$) would also be eigenvectors corresponding to $\lambda$.
\end{lemma}

\vspace{2mm}

\begin{proof}
	From definition of structural balance there exists $G_t$ such that
	\begin{align}
		G_tL_sG_t=L
		\label{eq:3}
		\quad\Rightarrow\quad
		\begin{cases}
			A_s^f=G'A^fG' \\
			B_s^f=\sigma_nG'B^f
		\end{cases}
	\end{align}
	where $G'=\mathrm{diag}(\sigma_1,\sigma_2,\dots,\sigma_{n-1})$.
	Without loss of generality, we can assume $\sigma_n=1$ since \eqref{eq:3} also holds for $-G$.
	This gives $B^f=G'B_s^f$ or $B^f=G'B_s^f$.

	From \cite{Rahmani2009a}, if $\mathcal{G}$ has the input symmetry structure, then there exists a nontrivial permutation matrix $J$ such that $JA^f=A^fJ$.
	
	Let $J'=G'JG'$, then the proofs of parts 1 and 2 follows:
	\begin{align*}
		J'A_s^f&=G'JG'G'A^fG'=G'JA^fG' \\
		&=G'A^fJG'=G'A^fG'G'JG'=A_s^fJ',
	\end{align*}
	and
	\begin{align*}
		J'^TB_s^f&=J'^TG'B^f=G'J^TG'G'B^f=G'J^TB^f\\
		&=-G'J^TA^f\textbf{1}=-G'A^fJ^T\textbf{1}=B_s^f.
	\end{align*}
	\noindent And for the last part of the proof
	\begin{align*}
		A_s^fv=\lambda v \quad\Rightarrow\quad
		A_s^fJ'v=J'A_s^fv=\lambda J'v,
	\end{align*}
	implying that $J'v$ is also an eigenvector for the same eigenvalue.
	Hence, $v-J'v$ would also be an eigenvector corresponding to $\lambda$.
\end{proof}

\vspace{2mm}

\begin{remark}
	\label{rmk:4}
	The matrix $J'$ introduced in Lemma 2 is in some sense correspondent to the permutation matrix $J$ in the unsigned case.
\end{remark}

\vspace{2mm}

\begin{theorem}
	\label{thm:2}
	The signed network system $\mathcal{G}_s$ is uncontrollable if it is input symmetric and structurally balanced.
\end{theorem}

\vspace{2mm}

\begin{proof}
	Let $(\lambda,v)$ be a pair of eigenvalue and eigenvector for $A_s^f$ so that $A_s^fv=\lambda v$. Then from Lemma \ref{lem:2}.3 we know that $v-J'v$ is also and eigenvector for $A_s^f$. Then
	\begin{align*}
		(v-J'v)^TB_s^f=v^TB_s^f-v^TJ'^TB_s^f=v^TB_s^f-v^TB_s^f=0,
	\end{align*}
	where both Lemma \ref{lem:2}.2 and PBH test are leveraged.
	The result implies that the system is uncontrollable.
\end{proof}

\vspace{2mm}

\begin{remark}
	\label{rmk:5}
	This result justifies viewing the gauge transformation as the unique invertible similarity transformation bridging state-space realizations $(A^f,B^f)$ and $(A_s^f,B_s^f)$, with $A^f=G'A_s^fG'$ and $B^f=G'B_s^f$.
	The focus of this paper, in the meantime, is to examine the graph-theoretic perspective of the signed networks and how structural balance paves the way for extending control theoretic analysis from unsigned to signed networks.
	This remark also applies to the results of the next sections.
\end{remark}

\vspace{2mm}

\begin{remark}
	\label{rmk:2}
	A signed symmetry implies the existence of an unsigned symmetry of $\mathcal{G}$. The converse is true when $\mathcal{G}_s$ is structurally balanced.
\end{remark}

\vspace{2mm}

\noindent\subsubsection{The MIMO Case}\label{sec:mathc-uncontr-if}

~

\vspace{1mm}

In this section, we examine how the notion of structural balance is interposed in the controllability analysis of multiple input signed networks. The results in this section are extensions to \cite{Rahmani2009a}. To this end, we leverage the machinery of interlacing and equitable partitions on graphs.

First, we modify two fundamental lemmas from~\cite{godsil2013algebraic} to the signed case and then provide the analysis which leads to sufficient conditions on the uncontrollability of the system.

\vspace{2mm}

\begin{definition}
	Let $G_t$ be the gauge transformation as in \eqref{eq:3}.
	Then $P'$ is the \emph{signed characteristic matrix} defined as $P'=G_tP$.
\end{definition}

\vspace{2mm}

\begin{lemma}
	\label{lem:3}
	Let $\pi$ be a parition of the structurally balanced signed graph $\mathcal{G}_s$, with adjacency matrix $A_s$ and signed characteristic matrix $P'$.
	Then $\pi$ is equitable if and only if the column space of $P'$ is $A_s$-\emph{invariant}.
\end{lemma}

\vspace{2mm}

\begin{proof} 
	\emph{(Necessity)} assume $\pi$ is equitable.
	Then from Lemma 9.3.1 in \cite{godsil2013algebraic}, $AP=P\hat{A}$ with $\hat{A}=A(\mathcal{G}_s/\pi)$.
	Thus, it follows from the structural balance: 
	\begin{align*}
		P\hat{A}=AP=G_tA_sG_tP
		\quad \Rightarrow \quad
		A_sP'=P'\hat{A}.
	\end{align*}
	\emph{(Sufficiency)} From Lemma 9.3.2 in \cite{godsil2013algebraic}, $\pi$ is equitable if there exists $B$ such that $AP=PB$.
	Then assuming that the column space of $P'$ is $A_s$-invariant, there exists $C$ such that
	\begin{align*}
		A_sP'=P'C
		\quad\Rightarrow\quad
		A_sG_tP=G_tPC
		\quad\Rightarrow\quad
		AP=PC
	\end{align*}
	Hence, $\pi$ is equitable.
\end{proof}

\vspace{2mm}

Lemma \ref{lem:3} shows how the gauge transformation is injected into the analysis of signed networks.

As discussed in \cite{Rahmani2009a}, we can now find an orthogonal decomposition of $\mathbb{R}^n$ using the signed characteristic matrix $P'$ as $\mathbb{R}^n=\mathcal{R}(P') \oplus \mathcal{R}(Q')$
where $\mathcal{R}(Q')=\mathcal{R}(P')^{\perp}$. Then an orthonormal basis for $\mathbb{R}^n$ can be formed as
\begin{align}
	T=[~\bar{P}'~|~\bar{Q}'~],
	\label{eq:4}
\end{align}
where $\bar{P}'$ and $\bar{Q}'$ represent normalized $P'$ and $Q'$ respectively and satisfy $\bar{P}'^T \bar{Q}'=0$, $\bar{P}'^T \bar{P}'=I$, and $\bar{Q}'^T \bar{Q}'=I$.

\vspace{2mm}

\begin{lemma}
	\label{lem:7}
	Given a connected signed graph $\mathcal{G}_s$, the system \eqref{eq:1} is uncontrollable if $L_s$ and $A_s^f$ share at least one common eigenvalue.
\end{lemma}

\vspace{2mm}

Lemma \ref{lem:7} is a derivation from Lemma 7.9 in \cite{Rahmani2009a}.
Since this is a general result depending on the Laplacian and its leading principal submatrix (floating graph), and the PBH test, the same holds for the signed case. 

From this point, the goal is to show that for some specific graph partition and the structural balance of the network, $L_s$ and $A_s^f$ share similar eigenvalues leading to uncontrollability.

\vspace{2mm} 

\begin{lemma}
	\label{lem:9}
	Suppose a structurally balanced signed graph $\mathcal{G}_s$ has an NEP $\pi$ with $\bar{P'}$ and $\bar{Q'}$ as in \eqref{eq:4}.
	Then the signed Laplacian $L_s$ is similar to the block diagonal matrix
	\begin{align*}
		\bar{L}_s=\begin{bmatrix}
		{L}_{P'}&\textbf{0}\\\textbf{0}&{L}_{Q'}
		\end{bmatrix},
	\end{align*}
	where ${L}_{P'}=\bar{P}'^TL_s\bar{P}'$ and ${L}_{Q'}=\bar{Q}'^TL_s\bar{Q}'$.
\end{lemma}

\vspace{2mm}

\begin{lemma}
	\label{lem:11}
	Let $\mathcal{G}_s^f$ be a signed floating graph, and $A_s^f$ be defined as in \eqref{eq:1} and $\bar{P'}$ and $\bar{Q'}$ be as in \eqref{eq:4}.
	If there exists an NEP $\pi_f$ in $\mathcal{G}_s^f$ and a $\pi$ in the original structurally balanced signed graph $\mathcal{G}_s$ such that all the nontrivial cells in $\pi_f$ are also cells in $\pi$, then $A_s^f$ is similar to the block diagonal matrix
	\begin{align*}
		\bar{A}_s^f=
		\begin{bmatrix}
		{A}_{P'}^f&\textbf{0}\\
		\textbf{0}&{A}_{Q'}^f
		\end{bmatrix},
	\end{align*}
	where ${A}_{P'}^f=\bar{P}_f'^TA_s^f\bar{P}_f'$ and ${A}_{Q'}^f=\bar{Q}_f'^TA_s^f\bar{Q}_f'$ with $\bar{P}_f'$ and $\bar{Q}_f'$ denoting the floating parts of $\bar{P}'$ and $\bar{Q}'$.
\end{lemma}

\vspace{2mm}

The proofs are similar to Lemmas 7.11, 7.12, and 7.14 in \cite{Rahmani2009a} and are skipped for succinctness.
One just needs to consider the role of structural balance and the fact that the signed characteristic matrix $P'$ needs to be replaced for $P$ due to the insertion of gauge transformation.

\vspace{2mm}

\begin{theorem}
	\label{thm:3}
	Given a connected structurally balanced signed graph $\mathcal{G}_s$ with the floating graph $\mathcal{G}_s^f$, the system in \eqref{eq:1} is uncontrollable if there exist NEPs on $\mathcal{G}_s$ and $\mathcal{G}_s^f$, $\pi$ and $\pi_f$, such that $\pi_f$ contains all nontrivial cells of $\pi$.
\end{theorem}

\vspace{2mm}

The main scheme of the proof is similar to Theorem 7.15 in \cite{Rahmani2009a}.
We modify the proof to show how the orthogonal basis formed by $\bar{P'}$ and $\bar{Q'}$ work in the new setup of signed networks. 

\vspace{2mm}

\begin{proof}
	Let $\bar{P'}$ and $\bar{Q'}$ be defined as in \eqref{eq:4}. 
	Following the convention in \cite{Rahmani2009a}, let $\pi\cap \pi_f=\{C_1,C_2,\dots,C_{r_1}\}$ with $|C_i|\geq2$ for $i=1,2,\dots,r_1$.
	Let the nontrivial cells contain the first $n_1$ nodes.
	Since $\pi_f$ contains all nontrivial cells of $\pi$, it follows
	\begin{align*}
		P'=\begin{bmatrix}
			P_1'&\textbf{0}\\
			\textbf{0}&I_{n-n_1}
		\end{bmatrix}_{n\times r}
		\quad \text{and} \quad
		P_f'=\begin{bmatrix}
			P_1'&\textbf{0}\\
			\textbf{0}&I_{n_f-n_1}
		\end{bmatrix}_{n_f\times r_f},
	\end{align*}
	where $P_1'\in\mathbb{R}^{n_1\times r_1}$ contains the nontrivial part of  the signed characteristic matrices.
	Let $\bar{P}'$ and $\bar{P}_f'$ be the normalization of $P'$ and $P_f'$ and define $\bar{Q}'$ and $\bar{Q}_f'$ as in \eqref{eq:4}.
	Then
	\begin{align*}
		\bar{Q}'=\begin{bmatrix}
			Q_1'\\\textbf{0}
		\end{bmatrix}_{n\times (n_1-r_1)}
		\quad , \quad
				\bar{Q}'_f=\begin{bmatrix}
				Q_1'\\\textbf{0}
				\end{bmatrix}_{n_f\times (n_1-r_1)},
	\end{align*}
	where $Q_1'\in\mathbb{R}^{n_1\times(n_1-r_1)}$ satisfies $Q_1'^TP_1=0$. It follows that $\bar{Q}_f'=R^T\bar{Q}'$ with $R=[I_{n_f},\textbf{0}]^T$.
	Then by Lemmas \ref{lem:9} and \ref{lem:11}, we get
	\begin{align*}
		\mathcal{L}_{Q'}=\bar{Q'}^TL_s\bar{Q'}=\bar{Q}_f'^TR^TL_sR\bar{Q}_f'=\bar{Q}_f'^TA_s^f\bar{Q}_f'=\mathcal{A}_{Q'}^f
	\end{align*}
	This implies that $L_s$ and $A_s^f$ share a block matrix and thus have at least one equal eigenvalue.
	Therefore, by Lemma \ref{lem:7} the system is uncontrollable.
\end{proof}

\vspace{2mm}

\begin{remark}
	\label{rmk:3}
	Theorem \ref{thm:3} provides a sufficient condition for uncontrollability of signed networks.
	However, this not a necessary condition; a system can be simultaneously uncontrollable and structurally unbalanced.
\end{remark}

\subsection{Stabilizability and Output Controllability}\label{sec:stab-outp-contr}

Recent developments in controllability have extended the idea of using symmetry to characterize controllability of linear systems.
Using Theorem~2 in \cite{Chapman2015a}, we can now characterize the controllability, output controllability, stabilizability and output stabilizability of the \emph{influenced signed consensus dynamics} mentioned in \eqref{eq:2}.

Our main result is the following theorem, which identifies structural balance as the key feature for uncontrollability on top of symmetry of the underlying unsigned graph.

\vspace{2mm}

\begin{theorem}
  \label{thr:2}
  Consider dynamics~\eqref{eq:2}. Let $J$ be a non-trivial signed fractional automorphism of $L$.
  Suppose further that $\mathcal{G}_s$ is structurally balanced with the signed Laplacian $L_s$ and gauge transformation $G_t$.
  Let $J_s=G_tJG_t$ with $B_s=G_tB(I)$ and $C_s=C(O)G_t$.
  Consider the following conditions: (a) $J_s B_s = B_s$, (b) $C_s(R)J_sC_s(\mathcal{V}\setminus R)^T = 0$,\\ $C_s(R)J_sC_s(R)^T=Z\neq I$, and (c) $J_s v_i = v_i$ for all $v_i\sim\lambda_i(L_s)>0$.
  Then
  \renewcommand{\theenumi}{\arabic{enumi}}
  \begin{enumerate}
    \item (a)  $\iff (-L_s,B)$ is uncontrollable
    \item (a) \& (b) $\iff (-L_s,B,C)$ is output uncontrollable
    \item (a) \& (c) $\iff (-L_s,B,C)$ is output unstabilizable
    \item (a) \& (b) \& (c) $\iff (-L_s,B,C)$ is output unstabilizable
  \end{enumerate}
  The results hold when $B\to B_s$ or $C\to C_s$.
\end{theorem}
\vspace{2mm}

The proof of Theorem~\ref{thr:2} is similar to the proof of Corollary~5 in \cite{Chapman2015a}, with several technical differences which we discuss here.
In particular, our result only requires a fractional automorphism $J$ of the underlying unsigned graph; the matrix $J_s$ does not need to be a signed fractional automorphism, which in this sense generalizes Corollary 5 of \cite{Chapman2015a}.

The following lemmas establish the equivalence of controllability under a gauge transformation of the $B$ and $C$ matrices for structurally balanced $L_s$, and some useful identities that will elucidate the role of the gauge transformation in Theorem~\ref{thr:2}.

\vspace{2mm}

\begin{lemma}
	\label{lem:4}
	Let $(L_s,B(I))$ be the pair in the dynamics~\eqref{eq:2},
	and let $B_s(I) = G_tB(I)$ for any gauge transformation $G_t$ (regardless of whether $L_s$ is structurally balanced).
	Then $(-L_s,B(I))$ is controllable if and only if $(-L_s,B_s(I))$ is controllable.
	Furthermore, letting $C_s(O) = C(O)G_t$, we have that $(-L_s,B(I),C(O))$ is output controllable if and only if $(-L,B_s(I),C_s(O))$ is output controllable.
\end{lemma}

\vspace{2mm}

\begin{lemma}
\label{lem:12}
Consider the dynamics in \eqref{eq:2}.
Suppose $L_s$ is structurally balanced with gauge $G_t$.
\begin{enumerate}
    \item Suppose that there is an automorphism $J$ such that $JL_s = L_sJ$.
      Then $J_s L_s = L_s J_s$, where $J_s = G_tJG_t$.
    \item Suppose $J$ is input symmetric. Then $J_sB_s = B_s$.
    \item Suppose that there exists $Z\neq I$ such that $ZC(O) = C(O)J$. Then $ZC_s(O) = C(O)J_s$.\label{lem3}
  \end{enumerate}
\end{lemma}

\vspace{1mm}

\begin{proof}
	\emph{(Lemma \ref{lem:4})} Note that the column space of the controllability matrix of $(-L_s,B)$ defined as $\mathcal{C}(B(I)):= \left[ 
	\begin{array}{cccc}
	B&-L_s B& \cdots& (-L_s)^{n-1}B
	\end{array}\right]$
	is spanned by columns of the form $(-L_s)^m e_i$ for $0\leq m\leq n-1$.         
	The action of a gauge $G_t$ on $B(I)$ is to multiply each column of $B(I)$ by $\pm 1$,
	and so the column space of the controllability matrix of $(-L_s,B_s)$ which is $\mathcal{C}(B_s(I)):= \left[ 
	\begin{array}{ccc}
	G_tB& \cdots& (-L_s)^{n-1}G_tB
	\end{array}\right]$
	is spanned by columns of the form $\sigma_i(-L_s)^m e_i$ for $0\leq m\leq n-1$, and $\sigma_i=\pm 1$.

	Clearly, $\mathrm{span}\{(-L_s)^m e_i\} = \mathrm{span}\{\sigma_i(-L_s)^m e_i\}$ and therefore $\mathrm{rank} [\mathcal{C}(B(I))] = \mathrm{rank} [\mathcal{C}(B_s(I))]$.
	The same argument applies for output controllability, considering the output controllability matrix $\mathcal{C}(B(I),C(O)) := \left[ 
	\begin{array}{ccc}
	CB&\cdots&C(-L_s)^{n-1}B 
	\end{array}\right]$.
\end{proof}

\vspace{2mm}

\begin{proof}
	\emph{(Lemma \ref{lem:12})} For each statement we have:
	
	\vspace{1mm}
	
	\noindent 1) $J_sL_s = G_tJG_t G_t L G_t = G_t JLG_t = G_t LJ G_t = G_tLG_tG_tJG_t = L_s J_s$

	\vspace{1mm}
	
	\noindent 2) $J_sB_s = G_sJG_sG_sB = G_sJB = G_sB=B_s$

	\vspace{1mm}
	
	\noindent 3) $ZC(O)G_t = C(O)JG_t = C(O)G_tG_tJG_t = C_s(O)J_s$.
\end{proof}

\vspace{1mm}

Using the equivalences established in these two lemmas, the proof of Theorem~\ref{thr:2} follows as the proof of Corollary 5 in~\cite{Chapman2015a}, but using $J_s$ instead of the signed fractional automorphism $P$.

\section{EXAMPLES}\label{sec:examples}

In this section, we show examples pertaining to the discussions in this paper.
In particular, we show the difference between consensus under signed and unsigned consensus.

\subsection{SISO Controllability}\label{sec:bipartite-consensus}

\begin{figure}[H]
 	\centering
 	\begin{minipage}{0.45\linewidth}
 		\centering
 		\subfloat[]{\label{fig:SB}\includegraphics[scale = 0.7]{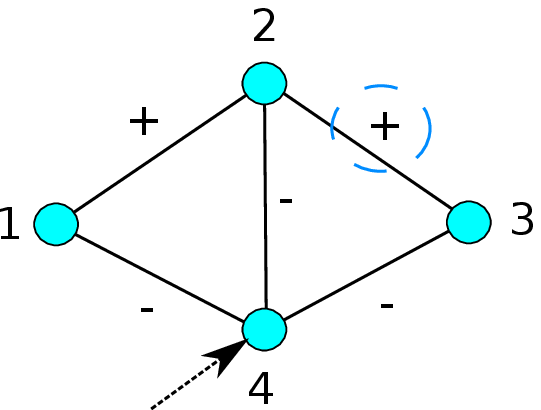}}
 	\end{minipage}
 	\hspace{3mm}
 	\begin{minipage}{0.45\linewidth}
 		\centering
 		\subfloat[]{\label{fig:SUB}\includegraphics[scale = 0.7]{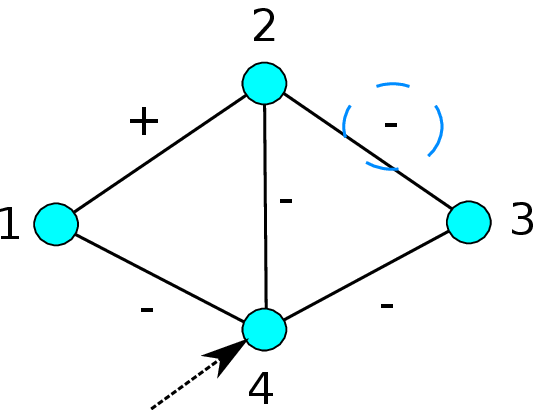}}
 	\end{minipage}
 	\caption{Signed graph with a single input symmetry about node 4. (a) Structurally balanced and uncontrollable (b) Structurally unbalanced and controllable}
 	\label{fig:lap}
\end{figure}

Consider the graphs in Figure~\ref{fig:lap}. Figure~\ref{fig:SB} shows the uncontrollability of the signed network as a consequence of structural balance and input symmetry. Figure~\ref{fig:SUB} verifies the role of structural balance on this result.
  
\subsection{MIMO Controllability}

In this example, we consider the MIMO case with two input signals injected onto nodes 4 and 5. The partition is equitable and $\pi=\{C_1,C_2,C_3,C_4\}$ and $\pi_f=\{C_1,C_2\}$.

\begin{figure}[H]
	\centering
	\begin{minipage}{0.45\linewidth}
		\centering
		\subfloat[]{\label{fig:MIMO1}\includegraphics[scale = 0.7]{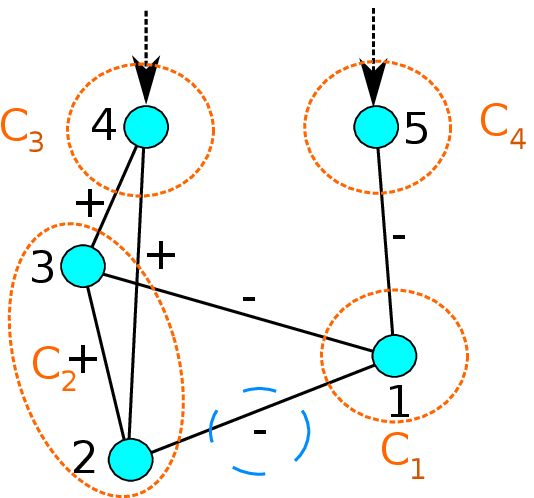}}
	\end{minipage}
	\hspace{3mm}
	\begin{minipage}{0.45\linewidth}
		\centering
		\subfloat[]{\label{fig:MIMO2}\includegraphics[scale = 0.7]{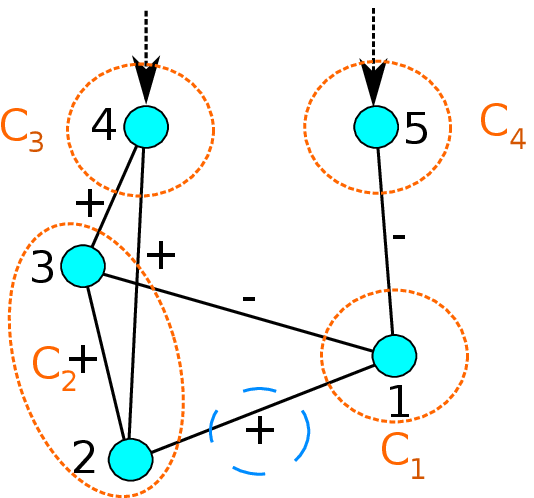}}
	\end{minipage}
	\caption{Signed graph with multiple input on nodes 4 and 5. (a) Structurally balanced and uncontrollable (b) Structurally unbalanced and controllable}
	\label{fig:MIMO}
\end{figure}

Then Figure~\ref{fig:MIMO} demonstrates how structural balance can influence the controllability of a signed network.

\subsection{Influenced Consensus}\label{sec:influenced-consensus}

Consider the influenced consensus of the networks in Figure~\ref{fig:lap}.
The network in Figure~\ref{fig:SB} is structurally balanced, but the one in Figure~\ref{fig:SUB} is not.

The controllability matrices for these two networks are, respectively:
\begin{align}
  \mathcal{C}_1 = \left[ 
  \begin{array}{cccc}
     0&     1 &    4 &   16\\
     0 &    1  &   4  &  16\\
     0  &   1   &  4   & 16\\
     1   &  3    &12    &48\\
  \end{array}\right],~\mathcal{C}_2 =  \left[ 
  \begin{array}{cccc}
     0 &    1 &    4 &   14\\
     0  &   1  &   6  &  32\\
     0   &  1   &  6   & 30\\
     1    & 3    &12    &52\\
  \end{array}\right].
\end{align}

As one can see, the network in Figure~\ref{fig:SB} is uncontrollable since $\mathcal{C}_1$ is rank-deficient, but the network in Figure~\ref{fig:SUB} is controllable, and hence one can conclude that unsigned symmetry is not sufficient for uncontrollability.

\section{CONCLUSIONS AND FUTURE WORKS}\label{sec:disc-concl}

In this paper, we characterized the controllability and stabilizability of signed consensus networks.
We showed that the topological notion of structural balance is a key construct pertaining to uncontrollability of signed networks.
In particular, structural balance induces a gauge transformation that permits the extension of the classical controllability and stabilizability analysis of consensus networks to signed networks.
We elucidated the role of structural balance both in the SISO and MIMO cases of the leader-follower (signed) consensus dynamics, and then extended these results to output controllability and stabilizability of the influenced networks.
Future works include extending the analysis to classes of nonlinear consensus networks and examining the application of controllability analysis to controller design and limits of performance for controlled networks.

\section{ACKNOWLEDGMENTS}
The authors would like to thank Dr. Airlie Chapman for enlightening conversations.

\bibliographystyle{IEEEtran}
\bibliography{ConsensusKoopman}

\end{document}